\documentclass[11pt]{amsart}
\usepackage{amsfonts}
\usepackage{cases}
\usepackage{mathrsfs}
\usepackage{bbm}
\usepackage{amssymb}
\usepackage{txfonts}
\usepackage{amscd}
\usepackage{amsfonts,latexsym,amsmath,amsthm,amsxtra,mathdots}
\usepackage[bookmarks,colorlinks]{hyperref}
\usepackage[all,cmtip]{xy}
\RequirePackage{amsmath} \RequirePackage{amssymb}
\usepackage{color}
\usepackage{colordvi}
\usepackage{multicol}
\usepackage{tikz}
\usepackage{hyperref}
\usepackage{graphicx}
\usepackage{amsmath}
\usepackage{amsmath,amscd}
\usetikzlibrary{matrix,arrows}

    \newcommand{\BQ}{{\mathbb {Q}}} \newcommand{\BR}{{\mathbb {R}}}

     \newcommand{\BZ}{{\mathbb {Z}}}

\def\-{^{-1}}

\newcommand{\delete}[1]{}

    \theoremstyle{plain}

\newtheorem{thm}{Theorem}[section]
\newtheorem{defn}[thm]{Definition} 
\newtheorem{lem}[thm]{Lemma}
\newtheorem{prop}[thm]{Proposition}
\newtheorem{cor}[thm]{Corollary}
\newtheorem{rem}[thm]{Remark}

    \numberwithin{equation}{section}

\def\Proof{\noindent{\bf Proof}\quad}
\def\qed{\hfill$\square$\smallskip}

\begin{document}

\title{The Farrell-Jones Conjecture for some nearly crystallographic groups}
\author{F. Thomas Farrell and Xiaolei Wu}

\begin{abstract}
In this paper, we
prove the K-theoretical and L-theoretical Farrell-Jones Conjecture with coefficients in an additive category for nearly crystallographic groups of the form $\BQ^n \rtimes \BZ$, where $\BZ$ acts on $\BQ^n$ as an irreducible integer matrix with determinant $d$, $|d |>1$.

\end{abstract}
\footnotetext{Date: August, 2014}



\maketitle


\section*{Introduction}

In the Farrell and  Linnell paper \cite{FL}, they proved that if the fibered isomorphism conjecture is true for all nearly crystallographic groups, then it is true for all virtually solvable groups. Recall that a nearly crystallographic group is a group of the form $A \rtimes C$ where $A$ is a torsion free abelian group of finite rank, $C$ is virtually cyclic and the action of $C$ on $A$ makes $A \otimes \BQ$ into an irreducible $\BQ C$ module (compare \cite{FL}, page 309).

Let $A  \in M(\BZ,n)$ with determinant $d$, $|d| > 1$, and $m(x)$ be its characteristic polynomial. Assume that $m(x)$ is irreducible over $\BQ$.  Let $K$ be $\mathbb{Q}(x) / {m(x)\mathbb{Q}(x)}$, then $K$ is a number field. Let ${\mathcal{O}}_K$ be its ring of integers. We define our group $\Gamma = {\mathcal{O}_K [\frac{1}{x}]}^+\rtimes_x \BZ$, where $\BZ$ acts on $\mathcal{O}_K [\frac{1}{x}]^+$ by multiplying $x$ (more details about this group will be given in Subsection \ref{anf}). In \cite{FW}, the authors proved the Farrell-Jones Conjecture with coefficients in an additive category  for all solvable Baumslag-Solitar groups. In the current paper we generalize the method there and prove the following (cf. Theorem \ref{main})

 \begin{flushleft} \textbf{Main Theorem}   The K- and L-theoretical Farrell-Jones Conjecture with coefficients in an additive category is true for the group $\Gamma$.
 \end{flushleft}

Note the truth of the Farrell-Jones Conjecture with coefficients in an additive category implies the fibered isomorphism conjecture.  For more information about the Farrell-Jones Conjecture and its fibered version, see for example \cite{FJ}. For precise formulation and discussion of Farrell-Jones Conjecture
with coefficients in an additive category, see for example \cite{BFL}, \cite{BLR2}.

As a corollary of our Main Theorem, we proved the following (cf. Corollary \ref{two})

 \begin{flushleft} \textbf{Corollary}
The K- and L-theoretical Farrell-Jones Conjecture with coefficients in an additive category is true for the group $\BQ^n \rtimes \BZ$, where $\BZ$ acts on $\BQ^n$ as an irreducible integer matrix with determinant $d$, $|d| > 1$.
\end{flushleft}

 \begin{flushleft} \textbf{Remark}
In the proof, we will stick to the case  $d > 1$ as it makes our notation cleaner.
\end{flushleft}

 \begin{flushleft} \textbf{Remark}  Independently, C. Wegner proved the Farrell-Jones Conjecture for all virtually solvable groups. The technique he uses is a combination of Farrell-Hsiang method and transfer reducibility while our method here does not use transfer reducibility. Also, our results does not depend on  Bartels and L\"uck's results on CAT(0) space \cite{BL2} \cite{BL3}. On the other hand, it seems our method here can not be applied to all nearly crystallographic groups, in particular, it can not be applied directly to the group $\BZ[\frac{1}{6}] \rtimes_{\frac{2}{3}} \BZ$.
\end{flushleft}

Our strategy is to show that $\Gamma$ is in fact a Farrell-Hsiang group, defined by  Bartels and L\"uck in \cite{BL1} . The main difficulty compared to our previous paper \cite{FW} is that we need to develop a theory for groups acting on trees. The second difficulty is that there is no canonical metric on the model $E(\Gamma)$ (cf. Section \ref{model}), so it is much more technical to prove some metric properties for $E(\Gamma)$.

The paper is organized as follows. In Section \ref{introfjc}, we give a brief introduction to the Farrell-Jones Conjecture, including  inheritance properties and the Farrell-Hsiang method. In Section \ref{gcotree}, we develop a theory of group actions on trees which can be used later. The construction  connects in natural to algebraic number theory. In Section \ref{model}, we produce a model $E(\Gamma)$ for the group $\Gamma$ and discuss some geometric properties of the space. In Section \ref{hfs}, we construct a horizontal flow space for the group $\Gamma$. Using results of Bartels, L\"uck and Reich, we produced a long-thin cover on this horizontal flow space. In the last section, we prove our main theorem.

In this paper, we will occasionally use FJC to abbreviate the K-theoretical or L-theoretical Farrell-Jones Conjecture with coefficients in an additive category. Let G be a (discrete) group acting on a space X. We say that the action is proper if for any $x \in X$ there is an open neighborhood U of $x$ such that $\{g \in G ~|~gU \bigcap U \neq {\emptyset}\}$ is finite.

\textbf{Acknowledgements.} This research was in part supported by the NSF grant DMS 1206622. The authors want to thank the Math Science Center at Tsinghua University for their warm hospitality during the preparation of this work.

\section{Some properties of the Farrell-Jones Conjecture}\label{introfjc}
In this section we give a brief introduction to the Farrell-Jones Conjecture.

We first list some useful inheritance properties valid for FJC. For details, see for example, \cite{BFL}, Section 1.3.
\begin{prop}\label{qut}
(1) If a group $G$
satisfies FJC, then every subgroup $H < G$ satisfies FJC.\\
(2) If $G_1$ and $G_2$ satisfy FJC, then the direct product $G_1 \oplus G_2$ and free product $G_1 \ast G_2$ satisfy FJC.\\
(3) Let $\{G_i ~|~ i \in I\}$ be a directed system of groups (with not necessarily injective
structure maps). If each $G_i$ satisfies FJC, then the direct limit $lim_{i \in I} G_i$
satisfies FJC.\\
(4) Let $\phi : G \rightarrow Q$ be a
group homomorphism. If $Q$ and $\phi^{-1}(C)$ satisfy FJC for every virtually  cyclic subgroup $C < Q$ then G satisfies FJC.\\
\end{prop}

\begin{lem}[Transitivity principle]\label{trp}
Let $\mathcal{F} \subseteq \mathcal{H}$ be two families of subgroups of G. Assume that G satisfies the K-Theoretic Farrell-Jones Conjecture with respect to $\mathcal{H}$ and that each $H \in \mathcal{H}$ satisfies the K-theoretic Farrell-Jones Conjecture with respect to $\mathcal{F}$. Then G satisfies the K-theoretic Farrell-Jones Conjecture with respect to $\mathcal{F}$. The same is true for the L-theoretic Farrell-Jones Conjecture.
\end{lem}

The following definition of Farrell-Hsiang group is taken from \cite{BL1}.
\begin{defn}\label{fhm}
Let $\mathcal{F}$ be a family of subgroups of the finitely generated group G. We call G a \textbf{Farrell-Hsiang Group} with respect to the family $\mathcal{F}$ if the following holds for a fixed word metric $d_G$:\\
  ~~ There exists a  fixed natural number $N$ such that for every natural number $m$ there is a surjective homomorphism $\Delta_m :G \rightarrow F_m$ with $F_m$ a finite group such that the following condition is satisfied. For every hyper-elementary subgroup H of $F_m$ we set $\bar{H} := \Delta^{-1}_m (H)$ and require that there exists a simplicial complex $E_H$ of dimension at most N with a cell preserving simplicial $\bar{H}$-action whose stabilizers belong to $\mathcal{F}$, and an $\bar{H}$-equivariant map $f_H: G \rightarrow E_H$ such that $d_G(g_0,g_1) < m$ implies $d^1_{E_H}(f_H(g_0),f_H(g_1)) < \frac{1}{m}$ for all $g_0,g_1 \in G$, where $d^1_{E_H}$ is the $l^1$-metric on $E_H$.

\end{defn}
\begin{rem}
As pointed out in \cite{BFL}, Remark 1.15, in order to check a group G is a Farrell-Hsiang group, it suffices to check these conditions for one hyper-elementary subgroup in every conjugacy class of such subgroups of $F_m$.
\end{rem}

Recall hyper-elementary group is defined as follows

\begin{defn}\label{he}
A \textbf{hyper-elementary group} H is an extension of a $p$-group by a cyclic group of order n, where $p$ is a prime number, $(n, p) = 1$, in other words,
there exists a short exact sequence
$$ 1 \rightarrow C_n \rightarrow H \rightarrow G_p \rightarrow 1 $$
where $C_n$ is a cyclic group of order n, $G_p$ is a $p$-group such that $(n, p) = 1$.
\end{defn}

With this definition, Bartels and L\"uck proved the following theorem:
\begin{thm}\label{fht}
Let G be a Farrell-Hsiang group with respect to the family $\mathcal{F}$. Then G satisfies the K-theoretic and L-theoretic Farrell-Jones Conjecture
with respect to the family $\mathcal{F}$.
\end{thm}

\section{Group acting on trees}\label{gcotree}

In this section, we develop a theory of group action on trees which can be used in the next section for constructing the model $E(\Gamma)$. It also plays a key role in the proof of our main theorem where we need  certain contracting maps to have an action on the tree. The basic construction is taken from Chapter II, \S 1 of Serre's Book \cite{Se}. Using Stallings folding and other techniques, we generalize results there to more general cases. The ideas already appeared in the appendix of \cite{FW}.

\subsection{Basic Construction}
We briefly recall the construction in Serre's book \cite{Se}, \S 1, Chapter II. Given a field\footnote{Here we require that multiplication in a field commutes while in Serre's book \cite{Se}, he does not require this.} $K$ with a discrete valuation $v$. Let $\mathcal{O}$ denote the valuation ring of $K$, i.e. the set of $x \in K$ such that $v(x) \geq 0$. We choose a uniformizer  $\pi$, i.e., an element $\pi \in K^\ast$ such that $v(\pi) = 1$ and let $k$ denote the residue field $\mathcal{O}/\pi\mathcal{O}$.

We let $V$ denote a vector space of dimension $2$ over $K$. A lattice of $V$ is any finitely generated $\mathcal{O}$-submodule of $V$ which generates the $K$-vector space $V$. If $x \in K^\ast$, and if $L$ is a lattice of $V$, $Lx$ is also a lattice of $V$. Thus the group $K^\ast$ acts on the set of lattices; we call the orbit of a lattice under this action its class. The set of lattice class is denoted by $T$.

Let $L$ and $L'$ be two lattices of V. By the invariant factor theorem there is an $\mathcal{O}$-basis $\{e_1,e_2\}$ of $L$ and integers $a,b$ such that $e_1 \pi^a, e_2 \pi^b$ is an $\mathcal{O}$-basis for $L'$. We define the distance between the corresponding lattice class $\Lambda$ and $\Lambda'$, $d(\Lambda, \Lambda') = |a - b|$. Two elements $\Lambda$, $\Lambda'$ of $X$ are called adjacent if $d(\Lambda, \Lambda') = 1$. In this way one defines a combinatorial graph structure on $T$. And the graph $T$ is a tree. Vertices of $T$ at distance $n$ from a fixed vertex $\Lambda_0$ correspond bijectively to direct factors of $L_0/L_0\pi^n$ of rank $1$, i.e. to points of the projective line $\textbf{P}(L_0/L_0 \pi^n) \cong \textbf{P}_1 ( \mathcal{O}/\pi^n \mathcal{O})$. In particular, edges with origin $\Lambda_0$ correspond bijectively to the points of $\textbf{P}(L_0/L_0 \pi) \cong \textbf{P}_1 ( k)$; if $q = Card(k)$, the number of these edges is $q+1$. For the purpose of this paper, $q$ will always be a finite number.

We let $GL(V)$ denote the group of $K$-automorphisms of $V$. We fix a basis $(e_1,e_2)$ for $V$, then $GL(V)$ is naturally isomorphic to $GL_2(K)$. Let $\mathfrak{p}$ be the prime ideal corresponding to the valuation function $v$.

\begin{defn}
We call a subgroup of $GL(V)$ $\mathfrak{p}$-bounded if it is a bounded subset of the vector space $End(V)$. If we identify $GL(V)$ with $GL_2(K)$, this means that there is an integer $l$ such that $v(s_{ij}) \geq l$ for each $s = (s_{ij}) \in G$.
\end{defn}

\subsection{Subgroups of $GL(V)$ and its action on the tree}  
\begin{figure}[h]
		\begin{center}
		\includegraphics[width=1.0\textwidth]{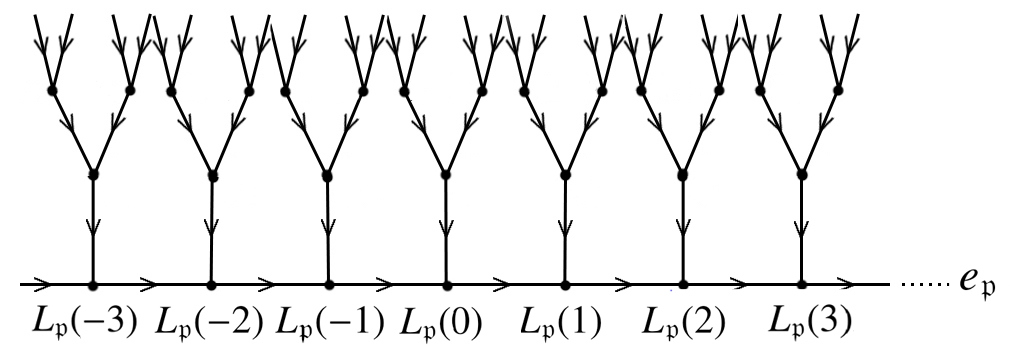}
		\end{center}
        \caption{The regular 3-valence tree $T_2$}
        \label{t2}
\end{figure}

We will use $v_{\mathfrak{p}}$ to denote the  valuation function corresponding  to  the prime ideal ${\mathfrak{p}}$. Choose a uniformizer $\pi_{\mathfrak{p}}$. We fix a basis $(e_1,e_2)$ for $V$. $GL(V)$ acts on the regular $(card(\mathcal{O}/\mathfrak{p}) +1)$-valence tree, which we will denote by $T_\mathfrak{p}$ (see Figure \ref{t2} for $T_2$).

 We are interested in the following subgroup of $GL(V)$ and its action on the tree.

$$G = \{ \left(
\begin{array}{cc}
    a & b\\
   0 & 1\\
    \end{array}\right) ~|~ a \in K^\ast, b\in K  \}$$

Note that $G$ is generated by elements of following two types: ${\left(
\begin{array}{cc}
    a  & 0 \\
   0 & 1 \\
    \end{array}\right)} $,
${\left(
\begin{array}{cc}
    1 & b \\
   0 & 1 \\
    \end{array}\right)} $. The action of $G$ on the tree $T_\mathfrak{p}$ has the following properties (compare, \cite{Se}, \S 1.3 Chapter II):
\begin{itemize}
\item
Let $L_{\mathfrak{p}}(n)$  represent the lattice class $ \mathcal{O} e_1 \oplus \pi_\mathfrak{p}^n \mathcal{O}e_2$, denote the end $L_{\mathfrak{p}}(0)L_{\mathfrak{p}}(1)L_{\mathfrak{p}}(2)\cdots$ as $e_\mathfrak{p}$, and the geodesic line $\cdots L_{\mathfrak{p}}(-2)L_{\mathfrak{p}}(-1) L_{\mathfrak{p}}(0)L_{\mathfrak{p}}(1)L_{\mathfrak{p}}(2)\cdots$ as $L_\mathfrak{p}$. Then the action of $G$ on $T_\mathfrak{p}$ fixes the end $e_\mathfrak{p}$. In particular, ${\left(
\begin{array}{cc}
    a  & 0 \\
   0 & 1 \\
    \end{array}\right)} $ acts on $L_{\mathfrak{p}}$ as translation by distance $-v_{\mathfrak{p}}(a)$. We can put an orientation on the tree $T_\mathfrak{p}$, where the orientation of every edge is pointing towards $e_\mathfrak{p}$ and $G$ preserves the orientation. Furthermore, we can define a Busemann function $f_\mathfrak{p}$ from $T_\mathfrak{p}$ to $\BR$ by
$$ f_{\mathfrak{p}}(P) := \lim_{n\rightarrow \infty} d(P, L_\mathfrak{p}(n)) -  n $$

 \item  At the vertex $L_{\mathfrak{p}}(0)$, there are $card(\mathcal{O}/\mathfrak{p})$ many vertices going into $L_{\mathfrak{p}}(0)$, which correspond bijectively  to elements in the residue field $\mathcal{O}/\mathfrak{p}$. For any $b \in \mathcal{O}$, ${\left(
\begin{array}{cc}
    1 & b \\
   0 & 1 \\
    \end{array}\right)} $ acts on these vertices by mapping $x$ to $x + b$, where $x \in \mathcal{O}/\mathfrak{p}$, and fixes $L_{\mathfrak{p}}(0)$, in fact it fixes the horoball $\{P \in T_{\mathfrak{p}}~|~ f_{\mathfrak{p}}(P) \leq 0\}$.

\item $G$ acts transitively on the vertices of $T_\mathfrak{p}$. The stabilizer of every vertex is bounded. The stabilizer of $L_{\mathfrak{p}}(0)$ consists of  matrices of the form
$ \left(
\begin{array}{cc}
    a & b\\
   0 & 1\\
    \end{array}\right)$, where $ v_{\mathfrak{p}}(a) = 0, v_{\mathfrak{p}}(b) \geq 0 $.

\end{itemize}

\begin{rem}\label{busm}
The Busemann function $f_d$ can be defined in the same way for any tree with a specified base point and a specified end.
\end{rem}

\subsection{Group action on $T_{\mathfrak{p}^k}$} \label{powertree}

Let $T_{\mathfrak{p}^k}$ be the regular $(card(\mathcal{O}/\mathfrak{p}^k) +1)$-valence tree.



 Let $G_{\mathfrak{p}^k}$ be the subgroup of $G$ generated by elements of the following three types:
 \begin{itemize}
 \item
type I: ${\left(
\begin{array}{cc}
    a  & 0 \\
   0 & 1 \\
    \end{array}\right)} $, $v_\mathfrak{p} (a) = k$;
    \item
   type II: ${\left(
\begin{array}{cc}
    m & 0 \\
   0 & 1 \\
    \end{array}\right)} $, $v_\mathfrak{p} (m) = 0$;
    \item
    type III: ${\left(
\begin{array}{cc}
    1 & b \\
   0 & 1 \\
    \end{array}\right)} $, $b \in \mathcal{O}$.
 \end{itemize}

 Matrices of type I act as translation on $T_\mathfrak{p}$ when restricted to the line $\cdots L_\mathfrak{p}(-2) \\ L_\mathfrak{p}(-1)L_\mathfrak{p}(0)L_\mathfrak{p}(1)L_\mathfrak{p}(2) \cdots$ by distance $-k$.   Recall $L_\mathfrak{p}(i) = \mathcal{O} e_1 \oplus \pi_\mathfrak{p}^i \mathcal{O}e_2$. Hence $f_\mathfrak{p}(M \circ P) = f_\mathfrak{p}(P) + k$,
for any point $P\in T_\mathfrak{p}$ and $M$ of type I matrix. When $M$ is matrix of type II and III, then $f_\mathfrak{p} (M \circ P) = f_\mathfrak{p} (P)$.

In order to get an action of $G_{\mathfrak{p}^k}$ on $T_{\mathfrak{p}^k}$ from its action on $T_{\mathfrak{p}}$, we need the following  definition of Stallings' folding, which is taken from \cite{BMF}.

\begin{defn}
Let T be an oriented tree with $G$ action.
Consider two edges $e_1$ and $e_2$ in T that are incident to a common vertex $v$. By
$\phi :e_1 \rightarrow e_2$ denote the linear homeomorphism fixing $v$. Then define an equivalence
relation " $\sim$ " on T as the smallest equivalence relation such that:\\
(i) $x \sim \phi(x)$, for all $x \in e_1$ , and\\
(ii) if $x \sim y$ and $g \in G$ then $g(x) \sim g(y)$.\\
The quotient space $T/\sim $ is a simplicial tree with a natural simplicial action of $G$.  Call the quotient map
$T \rightarrow T/\sim $ a fold.
\end{defn}

\begin{rem}
The key thing is that after the folding, G  has an induced action on the new tree.
\end{rem}

Consider the following subset of $T_\mathfrak{p}$
$$\{P \in T_\mathfrak{p} ~|~ -k +1 \leq f_{\mathfrak{p}}(P) \leq 0  \}$$
It has infinite many components. We will choose the component containing the vertex $L_\mathfrak{p}(0)$, which is a subtree of $T_\mathfrak{p}$, call it $\hat{T}_\mathfrak{p}$.

Proceeding similarly to the Stallings folding described in the appendix of  \cite{FW}, the Stalling folding on $T_\mathfrak{p}$ is generated by
folding every pair of edges in $\hat{T}_\mathfrak{p}$ that are incident to a common vertex in it. The resulting tree (compare Figure \ref{rep}) will be homeomorphic to $T_{\mathfrak{p}^k}$. Furthermore, if we delete all the $2$-valent vertices in the resulting tree, it will be exactly $T_{\mathfrak{p}^k}$. We will identify them as the same. $T_{\mathfrak{p}^k}$ will have a natural induced orientation, a specified geodesic line and a base point  under the identification. This further determines a Busemann function $f_{\mathfrak{p}^k}$ on $T_{\mathfrak{p}^k}$. We will denote the new geodesic line by $L_{\mathfrak{p}^k}$, and its vertices by $L_{\mathfrak{p}^k}(n)$. Now  $G_{\mathfrak{p}^k}$  has an induced action $T_{\mathfrak{p}^k}$ with the following properties:

\begin{figure}[h]
		\begin{center}
		\includegraphics[width=1.0\textwidth]{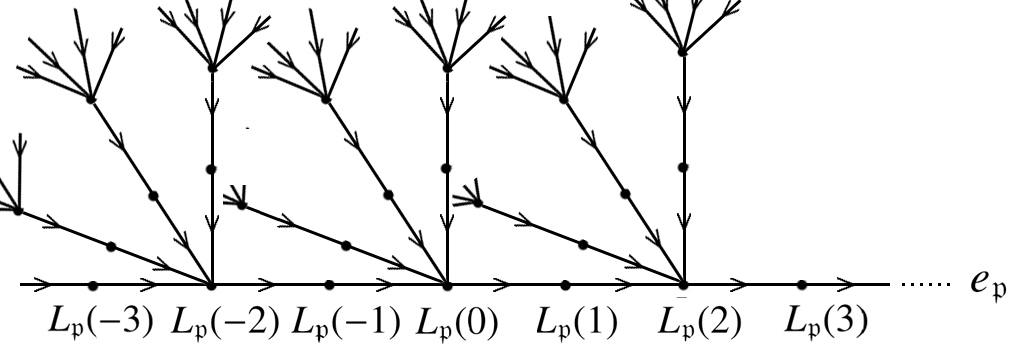}
		\end{center}
        \caption{$T_2$  after Stallings' folding}
        \label{rep}
\end{figure}

\begin{itemize}
 \item $G_{\mathfrak{p}^k}$ fixes the end $e_{\mathfrak{p}^k}$.

 \item Matrices of type I act as translation on $T_{\mathfrak{p}^k}$ by distance $-1$ when restricted to the geodesic line $L_{\mathfrak{p}^k}$. Hence $f_{\mathfrak{p}^k}(M \circ P) = f_\mathfrak{p}(P) + 1$,
for any point $P\in T_{\mathfrak{p}^k}$ and any matrix $M$ of type I. Matrices of type II fix the geodesic line $L_{\mathfrak{p}^k}$.  Furthermore, when $M$ is a matrix of type II or III, $f_{\mathfrak{p}^k} (M \circ P) = f_{\mathfrak{p}^k} (P)$.

  \item  At the vertex $L_{\mathfrak{p}^k}(0)$, there are $card(\mathcal{O}/\mathfrak{p}^k)$ many vertices going into $L_{\mathfrak{p}}(0)$, which correspond bijectively  to elements in the residue ring $\mathcal{O}/\mathfrak{p}^k$. For any $b \in \mathcal{O}$, ${\left(
\begin{array}{cc}
    1 & b \\
   0 & 1 \\
    \end{array}\right)} $ acts on these vertices by mapping $x$ to $x + b$, where $x \in \mathcal{O}/\mathfrak{p}^k$ and fixes $L_{\mathfrak{p}}(0)$,  in fact it fixes the horoball $\{P \in T_{\mathfrak{p}^k}~|~ f_{\mathfrak{p}^k}(P) \leq 0\}$.

  \item $G_{\mathfrak{p}^k}$ acts transitively on the vertices of the tree $T_{\mathfrak{p}^k}$. The stabilizer of every vertex is $\mathfrak{p}$-bounded.

\end{itemize}

\subsection{Action on the tree $T_{{{\mathfrak{p}_1}^{k_1}}{{\mathfrak{p}_2}^{k_2}}\ldots{{\mathfrak{p}_n}^{k_n}}}$} \label{treeaction}
In this subsection, we extend the group action on trees to more general cases. Let ${\mathfrak{p}_i}$ be the prime ideal corresponding to the valuation function $v_i$, and let $\mathcal{O}_i$ be the corresponding valuation ring. Now let $T_{{{\mathfrak{p}_1}^{k_1}}{{\mathfrak{p}_2}^{k_2}}\cdots{{\mathfrak{p}_n}^{k_n}}}$ be the $( \prod_{i=1}^{i=n}card(\mathcal{O}_i/{{\mathfrak{p}}_i}^{k_i}) +1)$-valence regular tree. And let $G_{{{\mathfrak{p}_1}^{k_1}}{{\mathfrak{p}_2}^{k_2}}\cdots{{\mathfrak{p}_n}^{k_n}}}$ be the subgroup of $G$ generated by elements of the following three types:
 \begin{itemize}
 \item
type I: ${\left(
\begin{array}{cc}
    a  & 0 \\
   0 & 1 \\
    \end{array}\right)} $, where $v_{\mathfrak{p}_i }(a) = k_i$, for every $1 \leq i \leq n$;
    \item
   type II: ${\left(
\begin{array}{cc}
    m & 0 \\
   0 & 1 \\
    \end{array}\right)} $, where $v_{\mathfrak{p}_i }(m) = 0$, for every $1 \leq i \leq n$;
    \item
    type III: ${\left(
\begin{array}{cc}
    1 & b \\
   0 & 1 \\
    \end{array}\right)} $, $b \in \bigcap_{i=1}^{i=n} \mathcal{O}_i$.
 \end{itemize}

Notice that
 $$G_{{{\mathfrak{p}_1}^{k_1}}{{\mathfrak{p}_2}^{k_2}}\cdots{{\mathfrak{p}_n}^{k_n}}} ~\subseteq~ G_{{\mathfrak{p}_1}^{k_1}}, ~for~ any~ 1 ~ \leq ~i ~\leq ~n$$

Hence as explained in subsection \ref{powertree}, $G_{{{\mathfrak{p}_1}^{k_1}}{{\mathfrak{p}_2}^{k_2}}\cdots{{\mathfrak{p}_n}^{k_n}}}$ acts on the tree $T_{{\mathfrak{p}_i}^{k_i}}$. Therefore it acts on the product space $ T_{{\mathfrak{p}_1}^{k_1}} \times T_{{\mathfrak{p}_2}^{k_2}} \times \cdots T_{{\mathfrak{p}_n}^{k_n
}}$ diagonally.  Now consider the ``diagonal subspace"
 $$\{(y_1,y_2,\cdots, y_{n})\in T_{{\mathfrak{p}_1}^{k_1}} \times T_{{\mathfrak{p}_2}^{k_2}} \times \cdots T_{{\mathfrak{p}_n}^{k_n
}}~ |~ f_{{{\mathfrak{p}_i}^{k_i}}}(y_i) = f_{{{\mathfrak{p}_j}^{k_j}}}(y_{j}),~ for ~any ~1\leq i , j\leq n \}$$
where $f_{{\mathfrak{p}_i}^{k_i}}$ is the Busemann function on the corresponding $T_{{\mathfrak{p}_i}^{k_i}}$. This "subspace" is left invariant under the $G_{{{\mathfrak{p}_1}^{k_1}}{{\mathfrak{p}_2}^{k_2}}\cdots{{\mathfrak{p}_n}^{k_n}}}$ action, hence it has an induced $G_{{{\mathfrak{p}_1}^{k_1}}{{\mathfrak{p}_2}^{k_2}}\cdots{{\mathfrak{p}_n}^{k_n}}}$-action.  It is not too hard to show that this subspace is homeomorphic to the tree $T_{{{\mathfrak{p}_1}^{k_1}}{{\mathfrak{p}_2}^{k_2}}\cdots{{\mathfrak{p}_n}^{k_n}}}$ we defined before, we will think of them as the same. Therefore we have an induced action of $G_{{{\mathfrak{p}_1}^{k_1}}{{\mathfrak{p}_2}^{k_2}}\cdots{{\mathfrak{p}_n}^{k_n}}}$ on $T_{{{\mathfrak{p}_1}^{k_1}}{{\mathfrak{p}_2}^{k_2}}\cdots{{\mathfrak{p}_n}^{k_n}}}$. There is a natural Busemann function defined on $T_{{{\mathfrak{p}_1}^{k_1}}{{\mathfrak{p}_2}^{k_2}}\cdots{{\mathfrak{p}_n}^{k_n}}}$  by $f_{_{{{\mathfrak{p}_1}^{k_1}}{{\mathfrak{p}_2}^{k_2}}\cdots{{\mathfrak{p}_n}^{k_n}}}}(y) = f_{{\mathfrak{p}_1}^{k_1}}(y_1)$. The geodesic line of $T_{{{\mathfrak{p}_1}^{k_1}}{{\mathfrak{p}_2}^{k_2}}\cdots{{\mathfrak{p}_n}^{k_n}}}$ corresponds to the ``diagonal subspace" of the product of the corresponding geodesic lines in $T_{{\mathfrak{p}_i}^{k_i}}$. We will represent this geodesic line by
$$\cdots L_{{{\mathfrak{p}_1}^{k_1}}{{\mathfrak{p}_2}^{k_2}}\cdots{{\mathfrak{p}_n}^{k_n}}}(-2)
L_{{{\mathfrak{p}_1}^{k_1}}{{\mathfrak{p}_2}^{k_2}}\cdots{{\mathfrak{p}_n}^{k_n}}}(-1)
L_{{{\mathfrak{p}_1}^{k_1}}{{\mathfrak{p}_2}^{k_2}}\cdots{{\mathfrak{p}_n}^{k_n}}}(0)
L_{{{\mathfrak{p}_1}^{k_1}}{{\mathfrak{p}_2}^{k_2}}\cdots{{\mathfrak{p}_n}^{k_n}}}(1)
L_{{{\mathfrak{p}_1}^{k_1}}{{\mathfrak{p}_2}^{k_2}}\cdots{{\mathfrak{p}_n}^{k_n}}}(2)
\cdots$$
For convenience, we will denote this geodesic line as $L_I$. There is a specified end of $T_{{{\mathfrak{p}_1}^{k_1}}{{\mathfrak{p}_2}^{k_2}}\cdots{{\mathfrak{p}_n}^{k_n}}}$,  $$L_{{{\mathfrak{p}_1}^{k_1}}{{\mathfrak{p}_2}^{k_2}}\cdots{{\mathfrak{p}_n}^{k_n}}}(0)
L_{{{\mathfrak{p}_1}^{k_1}}{{\mathfrak{p}_2}^{k_2}}\cdots{{\mathfrak{p}_n}^{k_n}}}(1)
L_{{{\mathfrak{p}_1}^{k_1}}{{\mathfrak{p}_2}^{k_2}}\cdots{{\mathfrak{p}_n}^{k_n}}}(2)
\cdots$$
which we will  denote as $e_I$.
The action of $G_{{{\mathfrak{p}_1}^{k_1}}{{\mathfrak{p}_2}^{k_2}}\cdots{{\mathfrak{p}_n}^{k_n}}}$ on $T_{{{\mathfrak{p}_1}^{k_1}}{{\mathfrak{p}_2}^{k_2}}\cdots{{\mathfrak{p}_n}^{k_n}}}$ has the following properties,

\begin{itemize}
 \item $G_{{{\mathfrak{p}_1}^{k_1}}{{\mathfrak{p}_2}^{k_2}}\cdots{{\mathfrak{p}_n}^{k_n}}}$ fixes the end $e_I$.

 \item Matrices of type I act as translation on $T_{{{\mathfrak{p}_1}^{k_1}}{{\mathfrak{p}_2}^{k_2}}\cdots{{\mathfrak{p}_n}^{k_n}}}$ by distance $-1$ when restricted to the geodesic line $L_I$. Matrices of type II fix the geodesic line $L_I$. Hence $f_{{{\mathfrak{p}_1}^{k_1}}{{\mathfrak{p}_2}^{k_2}}\cdots{{\mathfrak{p}_n}^{k_n}}}(M \circ P) = f_{{{\mathfrak{p}_1}^{k_1}}{{\mathfrak{p}_2}^{k_2}}\cdots{{\mathfrak{p}_n}^{k_n}}}(P) + 1$,
for any point $P\in T_{{{\mathfrak{p}_1}^{k_1}}{{\mathfrak{p}_2}^{k_2}}\cdots{{\mathfrak{p}_n}^{k_n}}}$ and any matrix $M$ of type I. If $M$ is a matrix of type II or III, then $f_{{{\mathfrak{p}_1}^{k_1}}{{\mathfrak{p}_2}^{k_2}}\cdots{{\mathfrak{p}_n}^{k_n}}}(M \circ P) = f_{{{\mathfrak{p}_1}^{k_1}}{{\mathfrak{p}_2}^{k_2}}\cdots{{\mathfrak{p}_n}^{k_n}}}(P)$.

  \item  At the vertex $L_{{{\mathfrak{p}_1}^{k_1}}{{\mathfrak{p}_2}^{k_2}}\cdots{{\mathfrak{p}_n}^{k_n}}}(0)$, there are $ \prod_{i=1}^{i=n}card(\mathcal{O}_i/{{\mathfrak{p}}_i}^{k_i})$ many vertices going towards it, which bijectively correspond to elements in the product of the residue ring $ \prod_{i=1}^{i =n}\mathcal{O}_i/{\mathfrak{p}_i}^k$. For any $b \in  \bigcap_{i=1}^{i=n} \mathcal{O}_i$, ${\left(
\begin{array}{cc}
    1 & b \\
   0 & 1 \\
    \end{array}\right)} $ acts on these vertices by mapping $(x_1,x_2,\cdots,x_n)$ to $(x_1 + b,x_2 + b,\cdots,x_n + b)$, where $x_i \in \mathcal{O}_i/{\mathfrak{p}_i}^k$ and fixes $L_{{{\mathfrak{p}_1}^{k_1}}{{\mathfrak{p}_2}^{k_2}}\cdots{{\mathfrak{p}_n}^{k_n}}}(0)$,  in fact it fixes the horoball $\{P \in T_{{{\mathfrak{p}_1}^{k_1}}{{\mathfrak{p}_2}^{k_2}}\cdots{{\mathfrak{p}_n}^{k_n}}}~|~ f_{{{\mathfrak{p}_1}^{k_1}}{{\mathfrak{p}_2}^{k_2}}\cdots{{\mathfrak{p}_n}^{k_n}}}(P) \leq 0\}$.

  \item $G_{{{\mathfrak{p}_1}^{k_1}}{{\mathfrak{p}_2}^{k_2}}\cdots{{\mathfrak{p}_n}^{k_n}}}$ acts transitively on the vertices of the tree $T_{{{\mathfrak{p}_1}^{k_1}}{{\mathfrak{p}_2}^{k_2}}\cdots{{\mathfrak{p}_n}^{k_n}}}$. And the stabilizer of each vertex is $({\mathfrak{p}_1}, {\mathfrak{p}_2}, \cdots{\mathfrak{p}_n})$-bounded. This means that there exists an integer $l$, such that $v_{\mathfrak{p}_i}(m_{st}) \geq l $ for each $M = (m_{st})$ in the stabilizer and for any $1 \leq i \leq n$.

\end{itemize}

\subsection{Applications to algebraic number fields} \label{anf}
Let $A  \in M(\BZ,n)$ with $det~A ~=~ d~ >~ 1$, and $m(x)$ be the characteristic polynomial of $A$. $m(x) = x^n + a_{n-1}x^{n-1} + \dots + a_1 x+ a_0$, $a_i \in \BZ$. Assume that $A$ is irreducible over $\BQ$. Then its characteristic polynomial equals its minimum polynomial.

Let $K$ be $\mathbb{Q}(x) / {m(x)\mathbb{Q}(x)}$, then $K$ is a number field. Let ${\mathcal{O}}_K$ be its ring of integers and $x{\mathcal{O}}_K$ be the principle ideal in ${\mathcal{O}}_K$ generated by $x$. One easily sees that $|{\mathcal{O}}_K/x{\mathcal{O}}_K| = d = |a_0|$.

Since ${\mathcal{O}}_K$ is a Dedekind domain, $x{\mathcal{O}}_K$ is a product of prime ideals, $\mathfrak{p}_1^{k_1} \mathfrak{p}_2^{k_2} \dots \mathfrak{p}_m^{k_m}$ where $k_i \in \BZ^+$ and  each $\mathfrak{p}_i$ is a prime ideal. Let $v_{\mathfrak{p}_i}$ be the valuation function corresponding to the prime ideal $\mathfrak{p}_i$. $v_{\mathfrak{p}_i}(x) = k_i$.

Let $G_x$ be the subgroup of $GL_2(K)$ generated by elements of the following three types:
 \begin{itemize}
 \item
type I: ${\left(
\begin{array}{cc}
    a  & 0 \\
   0 & 1 \\
    \end{array}\right)} $, where $v_{\mathfrak{p}_i }(a) = k_i$, for every $1 \leq i \leq n$;
    \item
   type II: ${\left(
\begin{array}{cc}
    m & 0 \\
   0 & 1 \\
    \end{array}\right)} $, where $v_{\mathfrak{p}_i }(m) = 0$, for every $1 \leq i \leq n$;
    \item
    type III: ${\left(
\begin{array}{cc}
    1 & b \\
   0 & 1 \\
    \end{array}\right)} $, $b \in \mathcal{O}_K$.
 \end{itemize}
As explained in the previous subsection, we have a $G_x$ action on a oriented regular (d+1)-valance tree with a specified Busemann function $f_x$, a specified geodesic line $\cdots P_x(-2)P_x(-1)P_x(0)P_x(1)P_x(2)\cdots$. We denote the geodesic line as $L_x$ and the end $P_x(0)P_x(1)P_x(2)\cdots$ as  $e_x$.

The action of $G_x$ has  the following properties:
\begin{itemize}
 \item $G_x$ fixes the end $e_x$.
 \item Matrices of type I act as translation on $T_d$ by distance $-1$ when restricted to the geodesic line $L_x$. Matrices of type II fixes the geodesic line $L_x$. Hence $f_x(M \circ P) = f_x(P) + 1$,
for any point $P\in T_d$ and any matrix $M$ of type I. When $M$ is a matrix of type II or III, then $f_x(M \circ P) = f_x(P)$.

  \item  At the vertex $P_x(0)$, there are $d$ vertices going into $P_x(0)$, which  bijectively correspond to elements in the residue ring $ \prod_{i=1}^{i =n}\mathcal{O}/{\mathfrak{p}_i}^k \cong \mathcal{O}/ \prod_{i=1}^{i =n}{\mathfrak{p}_i}^k \cong \mathcal{O}/x\mathcal{O}$. For any $b \in  \mathcal{O}$, ${\left(
\begin{array}{cc}
    1 & b \\
   0 & 1 \\
    \end{array}\right)} $   acts on these vertices by permutation, mapping $x$ to $x+ b$, where $x \in \mathcal{O}/x\mathcal{O}$ and fixes $P_x(0)$. In fact, it fixes the horoball $\{P \in T_{d}~|~ f_{x}(P) \leq 0\}$.

  \item $G_x$ acts transitively on the vertices of the tree $T_d$. And the stabilizer of each vertex is $({\mathfrak{p}_1}, {\mathfrak{p}_2}, \cdots{\mathfrak{p}_n})$-bounded. This means that there exists an integer $l$, such that $v_{\mathfrak{p}_i}(m_{st}) \geq l $ for each $M = (m_{st})$ in the stabilizer and for any $1 \leq i \leq n$.

\end{itemize}

Let $ \{~ \left(
\begin{array}{cc}
    x^k & b \\
   0 & 1 \\
    \end{array}\right)  ~|~ k \in \BZ, ~b \in {\mathcal{O}}_K[\frac{1}{x}]~ \}$ be a subgroup of $G_x$. Note that this subgroup is isomorphic to our group $\Gamma  = {\mathcal{O}_K [\frac{1}{x}]}^+\rtimes_x \BZ$. We have the following lemma.

\begin{lem} \label{actree}
 $\Gamma$ acts on the vertices of $T_d$ transitively and the stabilizer of every vertex is a free abelian group of rank $n$, where $n$ is the degree of the characteristic polynomial $m(x)$.
\end{lem}
\Proof The transitive part comes from the fact that $ \left(
\begin{array}{cc}
    x & 0 \\
   0 & 1 \\
    \end{array}\right)$  acts transitively on vertices of the geodesic line $L_x$ and the subgroup $\{~ \left(
\begin{array}{cc}
    1 & b \\
   0 & 1 \\
    \end{array}\right)  ~| ~b \in {\mathcal{O}}_K~ \}$ acts transitively on the $d$ edges going towards  the vertex $P_x(0)$. Since the action is transitive, we only need to show that the stabilizer of $P_x(0)$ is a free abelian group of rank $n$. Note first that the stabilizer of $P_x(0)$ is a subgroup of $\{~ \left(
\begin{array}{cc}
    1 & b \\
   0 & 1 \\
    \end{array}\right)  ~| ~b \in {\mathcal{O}}_K[\frac{1}{x}]~ \}$ and if $b \in {\mathcal{O}}_K$, then $ \left(
\begin{array}{cc}
    1 & b \\
   0 & 1 \\
    \end{array}\right)$  fixes $P_x(0)$. If $b \not\in {\mathcal{O}}_K$, then there exists a smallest $l > 0$ such that $b \in x^{-l}{\mathcal{O}}_K $, hence $b = x^{-l} a$ for some $a \in {\mathcal{O}}_K$, $a \not\in x{\mathcal{O}}_K$. Hence $\left(
\begin{array}{cc}
    1 & a \\
   0 & 1 \\
    \end{array}\right) $ fixes $P_x(0)$ and permutates (nontrivially) edges that point towards $P_x(0)$. On the other hand
    $$\left(
\begin{array}{cc}
    1 & x^{-l}a \\
   0 & 1 \\
    \end{array}\right) = \left(
\begin{array}{cc}
    x^{-l} & 0 \\
   0 & 1 \\
    \end{array}\right)\left(
\begin{array}{cc}
    1 & a \\
   0 & 1 \\
    \end{array}\right)\left(
\begin{array}{cc}
    x^l & 0 \\
   0 & 1 \\
    \end{array}\right)  $$
Since $\left(
\begin{array}{cc}
    x^l & 0 \\
   0 & 1 \\
    \end{array}\right)  $ acts on the geodesic line as translation by distance $-l$, $\left(
\begin{array}{cc}
    1 & x^{-l}a \\
   0 & 1 \\
    \end{array}\right)$ does not fix $P_x(0)$. Therefore we have proved the stabilizer of $P_x(0)$ is $\{~ \left(
\begin{array}{cc}
    1 & b \\
   0 & 1 \\
    \end{array}\right)  ~| ~b \in {\mathcal{O}}_K~ \} \cong {\mathcal{O}}_K^+$ which is a free abelian group of rank $n$.\qed

\section{A Model for $E(\Gamma)$} \label{model}

In this section, we construct a model for $E({\mathcal{O}}_K [\frac{1}{x}]^+ \rtimes_x \BZ)$, a contractible space
with free, proper and discontinuous $\Gamma$ action, where $\Gamma$ is the group $\mathcal{O}_K[\frac{1}{x}]^+ \rtimes_x \BZ$. We also put a metric on it, such that $\Gamma$ acts isometrically on
$E(\mathcal{O}_K[\frac{1}{x}]^+ \rtimes_x \BZ)$. Then we shall prove some inequalities that we will need in the future. Compare \cite{FW} for the case when $\Gamma$ is a solvable Baumslag-Solitar group.

Let $X = T_d \times \BR^n$, we claim there is a free properly discontinuous diagonal action of $\Gamma$ on it. The action of $\Gamma$ on $T_d$ comes from the fact that $\Gamma$ is a subgroup of $G_x$ (See Subsection \ref{anf}) via the following faithful representation
$$\varphi: \Gamma \rightarrow G_x \subset GL_2(K) ~ with ~ \varphi(b,k) = \left(
\begin{array}{cc}
    x^k & b \\
   0 & 1 \\
    \end{array}\right)$$
where $(b,k) \in \mathcal{O}_K[\frac{1}{x}]^+ \rtimes \BZ$. On the other hand, we also have an affine action of $\Gamma$ on $\BR^n$ given in the following way. Since $\mathcal{O}_K$ is a free abelian group of rank $n$, we can fix a basis for it, $e_1,e_2, \cdots e_n$. Since $x$ acts on $\mathcal{O}_K $ linearly, denote the corresponding matrix relative to the basis $e_1,e_2, \cdots e_n$ by $A_x$. For any $b \in \mathcal{O}_K[\frac{1}{x}]$, $b$ can be written uniquely as $b_1e_1 + b_2e_2 + \cdots + b_ne_n$, where $b_i \in \BZ[\frac{1}{d}]$ (since $det(A_x) = d$).

$$(b,k)(w_1,w_2, \cdots, w_n) = (w_1,w_2, \cdots, w_n) A^k_x + (b_1, b_2, \cdots, b_n)~~~~ \cdots~~~~(\ast)$$

\begin{rem}
Note here that  the actions of $\Gamma$ on $T_d$ and $\BR^n$ can  both be extended naturally to $G_x$; hence we have an induced diagonal action of $G_x$ on $T_d \times \BR^n$.
\end{rem}

\begin{prop}
$\Gamma$ acts on $T_d \times \BR^n$ properly, discontinuously and cocompactly.
\end{prop}
\Proof This is essentially a corollary of Lemma \ref{actree}. $\Gamma$ acts on the vertices of $T_d$ transitively while the stabilizer of $P_x(0)$ is a subgroup of $\mathcal{O}_K[\frac{1}{x}] \rtimes_x \{0\}$ and  isomorphic to $\BZ^n$. Note also that the action of $\mathcal{O}_K[\frac{1}{x}] \rtimes_x \{0\}$ on $\BR^n$ is simply by translation defined above by the formula ($\ast$), we have $\BR^n / \mathcal{O}_K[\frac{1}{x}] \rtimes_x \{0\}$ is homeomorphic to the $n$-dimensional torus $T^n$. Hence the action is proper, discontinuous and cocompact. \qed

One also sees that $T_d \times \BR^n /\Gamma $ is homeomorphic to a mapping torus of $T^n$. In fact, the linear action $A_x$ of $(x,0) \in \mathcal{O}_K[\frac{1}{x}] \rtimes_x  \BZ$ on $\BR^n$ induces an $d$-fold self covering map on the standard torus $T^n$ and $T_d \times \BR^n /\Gamma $ is homeomorphic to $T^n \times [0,1]/(w,0) \thicksim (w A_x,1)$. We need to put a metric on $T_d \times \BR^n$ so that $\Gamma$ acts on it by isometries. For that we only need to put a metric on the quotient space $T^n \times [0,1]/(w,0) \thicksim (w A_x,1)$. First put a standard flat metric on $T^n \times \{0\}$, and use the pullback metric by $A_x$ as the metric on $T^n \times \{1\}$. Then we linearly expand the metric to $T^n \times [0,1]$. Recall that $T_d$ has a standard metric with edge length $1$, $\BR^n$ has the standard Euclidean metric. Note that the metric we put here is different from the metric we constructed for the Baumslag-Solitar group in \cite{FW}, Section 2.



We are going to need the following lemmas in the future.
\begin{lem}{\label{dd}}
For any two points $(z_1,w_1),(z_2,w_2) \in T_d \times \mathbb{R}^n$,
$$d_{T_d \times \mathbb{R}^n}((z_1,w_1),(z_2,w_2)) \geq d_{T_d \times \mathbb{R}^n}((z_1,w_2),(z_2,w_2)) = d_{T_d}(z_1,z_2) $$
\end{lem}
\Proof The metric is not a warped product, but very close to one. In fact, any geodesic of the tree $T_d$ is naturally a geodesic of $T_d \times \BR^n$ (Compare Lemma 3.2 in \cite{CH} in the warped product case). Hence $d_{T_d \times \mathbb{R}^n}((z_1,w_2),(z_2,w_2)) = d_{T_d}(z_1,z_2)$.  The same argument for proving the warped product case in Lemma 3.1 in \cite{CH} shows $d_{T_d \times \mathbb{R}^n}((z_1,w_1),(z_2,w_2)) \geq d_{T_d}(z_1,z_2)$. \qed

\begin{rem}
It is not always true that $d_{T_d \times \mathbb{R}^n}((z_1,w_1),(z_2,w_2)) \geq d_{\BR^n}(w_1, w_2)$.
\end{rem}

\begin{cor}{\label{ddc}}
Let $z_1,z_2 \in T_d$, $w_1,w_2 \in \mathbb{R}^n$, then the following inequality holds
$$ d_{T_d \times \mathbb{R}^n}((z_1,w_1),(z_2,w_2)) \geq \frac{1}{2} d_{T_d \times \mathbb{R}^n}((z_1,w_1),(z_1,w_2))  $$
\end{cor}
\Proof  Combining the triangle inequality for metric spaces and  Lemma \ref{dd} yields the Corollary. \qed

\begin{lem}\label{rads}
Given a point $z_0 \in T_d$  then there exists an neighborhood $I$ of $z_0$ in $T_d$ containing $z_0$ such that for any $z \in I$ if $d_{\BR^n}(w_0,w_1)$ is greater than some fixed $\varepsilon_0$, then there exists $\varepsilon'_0$ such that $d_{T_d \times {\BR^n}}((z,w_0),(z,w_1)) > \varepsilon'_0$.
\end{lem}

\Proof The Lemma comes from the observation we can define another standard product metric in $T_d \times \BR^n$, the two metrics are different but they induce the same topology on $T_d \times \BR^n$. Hence they will induce the same topology on $\{z\} \times \BR^n$. In fact, the open ball of radius $\varepsilon_0$ centered at $(z,w_0)$ under the product metric is open under our metric also, hence there exists $\varepsilon'_0$ such that the ball of radius $\varepsilon'_0$  centered at $(z,w_0)$ under our metric is contained in the ball of radius $\varepsilon_0$ under the product metric. Therefore if $d_{\{z\} \times \BR^n}(w_0,w_1)$ is greater than $\varepsilon_0$, then $(z,w_1)$ is going to be outside the ball of radius $\varepsilon_0$ under the product metric centered at $(z,w_0)$. Hence it will be outside the ball of radius $\varepsilon'_0$ centered at $(z, w_0)$ under our metric, in particular  $d_{T_d \times {\BR^n}}((z,w_0),(z,w_1)) > \varepsilon'_0$. One can choose the neighborhood $I$ to be compact and sufficiently small, then we can further rechoose $\varepsilon'_0$ such that $d_{T_d \times {\BR^n}}((z,w_0),(z,w_1)) > \varepsilon'_0$ holds for any $z \in I$.\qed

\begin{lem}\label{vl}
Let $z_0$ be a fixed point in $T_d$, $w_1,w_2$ are two fixed points in $\mathbb{R}^n$, denote the distance $d_{T_d \times \mathbb{R}^n}((z_0, \frac{w_1}{n}),(z_0, \frac{w_2}{n})) $ by $D_n$, then for any $\epsilon>0$,  there exists an $N >0$ which depends on $D_1$ but is  independent of $z_0$ such that for any $n > N, D_n < \epsilon$. In particular,
 $$\lim_{n \rightarrow \infty} D_n = 0.$$
\end{lem}
\Proof  There are two metrics on $\{z_0\} \times \BR^n$ to measure the distance between $(z_0, \frac{w_1}{n})$ and $(z_0, \frac{w_2}{n})$, one is using geodesics in $T_d \times \BR^n$, the other one is using the standard Euclidean metric. Since they define the same topology on $\{z_0\} \times\BR^n$, we have for any ball of radius $\frac{\epsilon}{2} > 0$ centered at $(z_0,0)$ under our metric, we can find $\epsilon' > 0$ such that the ball of radius $\epsilon'$, under the Euclidean metric  is contained in the ball of radius $\frac{\epsilon}{2} $  under our metric. We first assume that $(z_0,w_1)$ lies in the fundamental domain (assume it contains $P_x(0)P_x(1) \times \{0\}$) which is compact, in particular $d_{ \mathbb{R}^n}(w_0, 0) < L$ for some $L>0$. Therefore both $(z_0, w_0)$ and $(z_0, w_1)$ lie in the ball of radius $D_1 + L$ centered at $(z_0, 0)$ under the Euclidean metric. Now we  choose $N$ big enough, such that for any $n> N$,  $(z_0, \frac{w_1}{n})$ and $(z_0, \frac{w_2}{n})$  lie in the ball of radius $\frac{\epsilon'}{2}$ under Euclidean metric, in particular their distance under the Euclidean metric is smaller than $\epsilon'$, hence both $(z_0,\frac{w_1}{n})$ and $(z_0,\frac{w_2}{n})$ will lie in the same ball of radius $\frac{\epsilon}{2}$ centered at $(z_0,0)$ under our metric, therefore the distance between $(z_0, \frac{w_1}{n})$ and $(z_0, \frac{w_2}{n})$ under our metric $D_n < \epsilon$. Since the fundamental domain is compact, We can rechoose $N$ so that for any $(z_0,w_1)$ lies in the fundamental domain the lemma holds.


Next, we extend the result to any $(z_0,w_1) \in T_d \times \BR^n$. Choose $g =(b,k) \in  \Gamma = {\mathcal{O}}_K [\frac{1}{x}]^+ \rtimes_x \BZ$ such that $g(z_0,w_1)$ lies in the fundamental domain. Then we have the following
$$d_{T_d \times \mathbb{R}^n}(g(z_0, \frac{w_1}{n}),g(z_0, \frac{w_2}{n})) = d_{T_d \times \mathbb{R}^n}((gz_0,  \frac{w_1}{n}A_x^k +b),(gz_0, \frac{w_2}{n}A_x^k + b))$$
Since translation in the $\BR^n$ direction does not change the distance,
$$d_{T_d \times \mathbb{R}^n}(g(z_0, \frac{w_1}{n}),g(z_0, \frac{w_2}{n})) = d_{T_d \times \mathbb{R}^n}((gz_0,  \frac{w_1}{n}A_x^k + \frac{b}{n}),(gz_0, \frac{w_2}{n}A_x^k + \frac{b}{n})) $$
$$\hspace{29.5mm} = d_{T_d \times \mathbb{R}^n}((gz_0,  \frac{1}{n} gw_1),(gz_0, \frac{1}{n}gw_2))$$

In particular, when $n=1$, $d_{T_d \times \mathbb{R}^n}((gz_0,  \frac{1}{n} gw_1),(gz_0, \frac{1}{n}gw_2)) = D_1$ and $(gz_0,  \frac{1}{n} gw_1) $ lies in the fundamental domain, therefore when $n >N$,  \\$d_{T_d \times \mathbb{R}^n}((gz_0,  \frac{1}{n} gw_1),(gz_0, \frac{1}{n}gw_2)) < \epsilon$. Hence $d_{T_d \times \mathbb{R}^n}(g(z_0, \frac{w_1}{n}),g(z_0, \frac{w_2}{n})) <  \epsilon$. $g$ acts on $T_d \times \BR^n$ as isometry, hence $d_{T_d \times \mathbb{R}^n}((z_0, \frac{w_1}{n}),(z_0, \frac{w_2}{n})) <  \epsilon$. This finishes our proof.  \qed

\begin{lem}\label{control}
Let $T$ be a positive integer, $z_0, z_1 \in T_d$, $z_1$ lies in the geodesic ray connecting $z_0$ and the end $e_x$, $d_{T_d} (z_0, z_1) = T$. $z$ lies in the geodesic connecting $z_0$ and $z_1$. Then there exists $\beta > 1$ depends on $T$ and independent of $z_1$, such that for any $w\in \BR^n$, if $d_{T_d \times \BR^n} ((z_1,w),(z_1,0)) = R \leq 1$, then $d_{T_d \times \BR^n} ((z,w),(z,0)) \leq \beta R$.
\end{lem}
\Proof One first fix $w$ and $z_1$, then let $M$ be the maximum value of $d_{T_d \times \BR^n} ((z,w),(z,0))$ where $z$ lies in the geodesic connecting $z_0$ and $z_1$ (which is compact). Then we can just choose $\beta = \frac{M}{R}$ in this case (assume $R>0$). Moreover we can choose such $\beta$ for all $w$ such that $d_{T_d \times \BR^n} ((z_1,w),(z_1,0)) \leq 1$ at the same time since the set of all such $w$ forms a compact set. Notice that when $w$ getting close to $0$, the metric getting close to Euclidean metric and we are able to choose such $\beta$ for those  $w$ continuously (one can also  extend $\beta$ to those $w$ by simply making $\beta$ large enough). One can further assume $z_1$ to lie on the edge $P_x(0)P_x(1)$ by using isometry  action of $\Gamma$  on $T_d\times \BR^n$  to translate any $z_1$  to $P_x(0)P_x(1)$. In fact we can always translate the geodesic connecting $z_0$ and $z_1$ to a geodesic lies in the specified geodesic line $L_x$ and $z_1$ lies in $P_x(0)P_x(1)$. Finally, because  $P_x(0)P_x(1)$ is compact, we can further choose a $\beta$ works for all $z_1 \in P_x(0)P_x(1)$. Hence $\beta$ is independent of $z_1$.

\section{Flow Space for $ T_d \times {\mathbb{R}^n}$}\label{hfs}
In this section we define a flow space for $E(\mathcal{O}_K[\frac{1}{x}] \rtimes \BZ)$. The construction here is parallel to the one in our previous paper \cite{FW}, but as the metric we put on $T_d \times \BR^n$ is not as nice as before, so we need some new ideas here in certain proofs.

We first introduce Bartels and L\"uck's flow space starting with the notion of generalized geodesic.
\begin{defn}
Let X be a metric space. A continuous map $c : \mathbb{R} \rightarrow X$ is called a \textbf{generalized geodesic} if
there are $c_{-}, c_{+} \in \bar{\mathbb{R}} :=  \mathbb{R} \coprod \{-\infty, \infty \}$ satisfying

$$  c_- \leq c_+, c_- \neq \infty, c_+ \neq -\infty  $$
such that $c$ is locally constant on the complement of the interval $I_c := (c_-, c_+)$ and
restricts to an isometry on $I_c$.
\end {defn}

\begin{defn}
Let $(X, d_X)$ be a metric space . Let \textbf{$FS(X)$} be the set of all generalized geodesics in
$X$. We define a metric on $FS(X)$ by

$$    d_{FS(X)}(c,d) := \int_{\mathbb{R}} \frac{d_X(c(t), d(t))}{2e^{|t|}} dt             $$
\end{defn}

The flow on $FS(X)$ is defined by
$$ \Phi: FS(X) \times \mathbb{R} \rightarrow FS(X) $$
where $ {\Phi}_{\tau} (c)(t) = c(t + \tau)$ for $ \tau \in \mathbb{R}$, $c \in FS(X)$ and $t \in \mathbb{R}$.

\begin{lem}{\label{inq}}
The map $\Phi$ is a continuous flow and if we let $c,d \in FS(X)$, $\tau \in \mathbb{R}$,
then the following inequality holds
$$ e^{-|\tau|} d_{FS(X)}(c,d) \leq d_{FS(X)}(\Phi_{\tau}(c),\Phi_{\tau}(d)) \leq e^{|\tau|} d_{FS(X)}(c,d)  $$
\end{lem}
\Proof A more general version is proved  in  \cite{BL2},  Lemma 1.3. \qed

Note that the isometry group of $(X, d_X)$ acts canonically on $FS(X)$. Recall a  map is proper if the inverse image of every compact subset is compact. Bartels and L\"uck  also proved the following for the flow space $FS(X)$ in \cite{BL2} Proposition 1.9 and 1.11.

\begin{prop} \label{ac}
If $(X,d_X)$ is a proper metric space, then $(FS(X),d_{FS(X)})$ is a proper metric space, in particular it is a complete metric space. Furthermore, if a group $\Gamma$ acts isometrically and properly on $(X, d_X)$, then $\Gamma$ also acts on $(FS(X),d_{FS(X)})$ isometrically and properly. In addition, if $\Gamma$ acts cocompactly on $X$, then $\Gamma$ acts cocompactly on $FS(X)$.

\end{prop}

Now we define our flow space by
 $$HFS(T_d \times\mathbb {R}^n) := FS(T_d)\times \mathbb {R}^n$$
  where $T_d$ has its natural metric with edge length $1$.   Since $\Gamma$ has an action on both $FS(T_d)$ and $\mathbb{R}^n$,
 $\Gamma$ will have a diagonal action on $FS(T_d)\times \mathbb{R}^n$ also. One can think of $HFS(T_d \times\mathbb {R}^n)$ as the horizontal subspace of $FS(T_d \times\mathbb {R}^n)$. In fact, there is a natural  embedding of  $HFS(T_d \times\mathbb {R}^n)$ (as a topological space with product topology) into $FS(T_d \times \mathbb {R}^n)$ defined as follows:
for a generalized geodesic $c$ on $T_d$, and $w \in \mathbb{R}^n$, we define a generalized geodesic on $T_d \times \mathbb{R}^n$, which maps $t \in \mathbb{R}$ to $(c(t),w) \in T_d \times \mathbb {R}^n$ (by Lemma \ref{dd}, for $t \in [0,\infty]$, ((c(t),w)) is a geodesic of $T_d\times \BR^n$). $HFS(T_d \times\mathbb {R}^n)$ will inherit a metric from this embedding.

For the rest of this section, let $X = T_d \times \mathbb {R}$.

\begin{lem}\label{prop}
The flow space $HFS(T_d \times\mathbb {R}^n)$ is a proper metric space, in particular a complete metric space.
\end{lem}

\Proof The proof here is parallel to the proof of Lemma 3.5 in \cite{FW} where in the proof a key lemma (Lemma 2.6) there is replaced by Lemma \ref{rads}. For completeness, we write down the proof here.
In order to prove $HFS(T_d \times\mathbb {R}^n)$ is a proper metric space, we need to show every closed ball $B_r(c) ~= ~\{c' ~|~ d_{HFS(T_d \times\mathbb {R}^n)}(c,c') \leq r\}$ in $HFS(T_d \times\mathbb {R}^n)$ is compact. Let $\{c_i\}$ be a Cauchy sequence in the closed ball $B_r(c)$, we need to show it converges to a point in $B_r(c)$. Since the space $FS(T_d \times \mathbb {R}^n)$ is proper, we can now assume $\{c_i\}$ converges to a point $c_0$ in $FS(T_d \times \mathbb {R}^n)$. We only need to show $c_0 \in HFS(T_d \times\mathbb {R}^n)$. Denote the projection map from $T_d \times \mathbb{R}^n$ to $T_d$ as $q_1$, from $T_d \times \mathbb{R}^n$ to $\mathbb{R}^n$ as $q_2$, then $c_i(t) = (q_1(c_i(t),q_2(c_i(t)))$. Suppose $c_0 \notin HFS(T_d \times\mathbb {R}^n)$, then $q_2(c_0(t))$ is not a constant map. Choose a big enough close interval $I$ in $\mathbb{R}$ such that $q_2(c_0(t))$ restricted to $I$ is not a constant. In particular, we can choose two point $t_1, t_2 \in I$, such that $q_2(c(t_1)) \neq q_2(c(t_2))$. Let $A_1 = q_2(c(t_1))$ and $A_2 = q_2(c(t_2))$. Let $\delta = |A_1 - A_2|$ which is  the distance between $A_1$ and $A_2$ under the standard Euclidean metric. Let $I_1 = \{t \in I~|~ |q_2(c_0(t)) - A_1| \leq \frac{\delta}{4}\}$, correspondingly $I_2 = \{t \in I~|~ |q_2(c_0(t)) - A_2| \leq \frac{\delta}{4}\}$. Note $I_1$ and $I_2$ are nonempty compact sets with measure bigger than $0$. Now for any given $c_i$, if $|q_2(c_i) - A_2|\geq \frac{\delta}{2}$, then for $t \in I_2$, $d_{\BR^n}(q_2(c_0(t)), q_2(c_i(t))) \geq \frac{\delta}{4}$ and

$$d_{FS(X)}(c_0, c_i) = \int_{\mathbb{R}} \frac{d_{T_d \times \mathbb{R}^n}(c_0(t), c_i(t))}{2e^{|t|}} dt \hspace*{69mm}$$

$$\geq \int_{I_2} \frac{d_{T_d \times \mathbb{R}^n}(c_0(t), c_i(t))}{2e^{|t|}} dt \hspace*{46mm}$$
$$\geq \int_{I_2} \frac{d_{T_d \times \mathbb{R}^n}((q_1(c_0(t)),q_2(c_0(t))), (q_1(c_i(t)),q_2(c_i(t))))}{2e^{|t|}} dt $$
$$\hspace*{13mm} = \int_{I_2} \frac{\frac{1}{2}d_{T_d \times \mathbb{R}^n}(c_0(t), (q_1(c_0(t)), q_2(c_i(t))))}{2e^{|t|}} dt\hspace*{5mm}( by~ Corollary~ \ref{ddc})$$

Note that here $d_{\BR^n}(q_2(c_0(t)), q_2(c_i(t))) \geq \frac{\delta}{4}$, we can choose $I'_2 \subset I_2$ such that Lemma \ref{rads} holds, hence there exists $\delta'$ such that $d_{T_d \times \mathbb{R}^n}(c_0(t), (q_1(c_0(t)), q_2(c_i(t)))) > \delta'$ for any $t \in I'_2$ then
$$  \int_{I_2} \frac{\frac{1}{2}d_{T_d \times \mathbb{R}^n}(c_0(t), (q_1(c_0(t)), q_2(c_i(t))))}{2e^{|t|}} dt ~~~ \geq ~~~ \int_{I'_2} \frac{\frac{1}{2}\delta'}{2e^{|t|}} dt > 0$$

The last integral is independent of $c_i$; denote its value as $\epsilon_1$. If $|q_2(c_i) - A_2 | \leq \frac{\delta}{2}$, then $|q_2(c_i) - A_1|\geq \frac{\delta}{2}$. For the same reason,  there exists  $\epsilon_2 > 0$ such that $d_{FS(X)}(c_0, c_i) \geq \epsilon_2$. Let $\epsilon = \min(\epsilon_1,\epsilon_2) > 0$, then $d_{FS(X)}(c_0, c_i) \geq \epsilon >0 $.
Hence the sequence $\{c_i\}$ can never converge to $c_0$, contradiction. \qed

\begin{rem}
The proof in fact shows that the embedding $HFS(X) \subset FS(X)$ is a closed $\Gamma$-equivariant embedding.
\end{rem}

We define now the flow
$$ \Phi: HFS(T_d \times\mathbb{R}^n )\times \mathbb{R} \rightarrow HFS(T_d \times\mathbb{R}^n) $$
by $ {\Phi}_{\tau} ((c,w))(t) = (c(t + \tau),w)$ for $c \in FS(T_d)$ and $ \tau,w,t \in \mathbb{R}$. Note $\Phi$ is a $\Gamma$-equivariant flow.
~\\

\begin{lem} {\label{l3}}
The flow space $HFS(X)$ has the following properties:
\begin{itemize}
    \item[(i)] $\Gamma$ acts properly and cocompactly on $HFS(X)$.
    \item [(ii)]  Given $C > 0$, there are only finitely many $\Gamma$ orbits of periodic flow curves with period less than $C$ (but bigger than $0$).
    \item [(iii)] Let $HFS(X)^{\mathbb{R}}$ denote the $\mathbb{R}$-fixed point set, i.e., the set of points $c \in HFS(X)$ for which ${\Phi}_{\tau} (c) = c$ for all $\tau \in \mathbb{R}$, then $HFS(X) - HFS(X)^{\mathbb{R}}$ is locally connected.
    \item[(iv)] If we put
    $$ k_\Gamma := sup\{ |H| | H \subseteq \Gamma subgroup~with ~finite~ order~ |H| \} ;$$
    $$d_{HFS(X)} := dim(HFS(X) - {HFS(X)}^{\mathbb{R}}); \hspace*{37mm}$$
    then $k_\Gamma$ and $d_{HFS(X)}$ are finite.

\end{itemize}

\end{lem}

\Proof We skip the proof of (i), (iii) and (iv) here  as they are parallel to the proof of Lemma 3.7 of \cite{FW}. The proof of (ii) needs some modification and uses the fact that $x$ acts on $\mathcal{O}_K[\frac{1}{x}]$ as an irreducible matrix (over $\BQ$) with determinant greater than $1$.

 Proof of (ii). Note that periodic orbits in $HFS(T_d\times \mathbb{R}^n)$ are periodic orbits in $FS(T_d\times \mathbb{R}^n)$, which move horizontally (i.e., move along the tree direction, with $\mathbb{R}^n$ coordinate fixed). Note also that the embedding of $HFS(T_d\times \mathbb{R}^n)$ into $FS(T_d\times \mathbb{R}^n)$ is a $\Gamma$-equivariant map, and there are only finitely many nonzero horizontal periodic geodesics on $(T_d\times \mathbb{R}^n)/{\Gamma}$ of period less than $C$. In fact $(T_d\times\mathbb{R}^n)/{\Gamma} = T^n \times [0,1]/(w,0) \thicksim ( A_x w,1)$ where $T^n$ is the n-dimensional torus (see Section \ref{model} for more details) and horizontal periodic geodesics with period $m$ on it correspond to the number of solutions of the equation $A_x^m w ~\equiv~ w~ (~mod~ (1,1,\cdots,1)~)$ times $d^m$ (number of branchings in $T_d$), where $w \in (\mathbb{R}/\mathbb{Z})^n$. Since for any positive integer $m$, the equation has finitely many solutions. In fact the solutions are $w = (A_x^m -I)^{-1} k, ~k \in \mathbb{Z}^n$, $w\in (\mathbb{R}/\mathbb{Z})^n$. Here $A_x^m -I$ is invertible due to the fact that $A_x$ is an irreducible $k \times k$ matrix (over $\BQ$) with $k >1$.  So there will be only finitely many nonzero horizontal periodical orbits on $HFS((T_d\times\mathbb{R})/{\Gamma} )$ with period less than $C$. Hence the claim in (ii) now follows.

\qed

\begin{rem} \label{remfl}
We define an embedding $\Psi: T_d \times \mathbb{R}^n \rightarrow FS(T_d) \times \mathbb{R}^n$, by $(z,w) \rightarrow (c_z,w)$, where $c_z$
is the unique generalized geodesic which sends $(-\infty, 0)$ to $z$, and $[0,\infty)$ isometrically to the geodesic $[z,e_x)$ where $e_x$ is the specified end of $T_d$. Also, we can flow this embedding by flowing its image in $HFS(X)$; define $\Psi_{\tau}(z,w) = \Phi_{\tau}(\Psi(z,w))$. It is easy to see that $\Psi_{\tau}$ is a $\Gamma$-equivariant map since $e_x$ is fixed under the group action.
\end{rem}

We need the following lemma in the proof of our main theorem.
\begin{lem} \label{lpar}
Let $z_0$ be a fixed point in $T_d$, $w_1,w_2$ are two fixed points in $\mathbb{R}^n$, and $P_n = (z_0,\frac{w_1}{n})$, $Q_n = (z_0,\frac{w_2}{n})$, $d_X(P_1,Q_1) < D$. Then for any $\epsilon > 0$, there exists a number $\bar{N}$, which depends only on $\epsilon$, $D$ and $d$, such that for any $n > \bar{N}$
$$d_X(P_n,Q_n) < \frac{\epsilon}{4}$$
and
$$ d_{HFS(X)}(\Psi(P_n),\Psi(Q_n)) \leq \epsilon$$
\end{lem}
\Proof Choose $T$ to be the first positive integer greater than $\ln {\frac{4}{\epsilon}}$. Since $\Psi(P_n)(T)$ and $\Psi(Q_n)(T)$ have the same $T_d$ coordinate, by Lemma \ref{vl}, we can choose a big enough integer $\bar{N}$ such that for any $n > \bar{N}$, $d_X(\Psi(P_n)(T),\Psi(Q_n)(T)) < \frac{\epsilon}{4\beta}$, where $\beta >1$ is determined by Lemma \ref{control} for our $T$ here. Note $\bar{N}$ depends only on $\epsilon$, $D$, $T$ (which is determined by $\epsilon$) and $d$.  Using the definition of generalized geodesic, we have $d_{X}( \Psi(P_n)(t), \Psi(P_n)(T)) = t-T$ and $d_{X}(\Psi(Q_n)(t), \Psi(Q_n)(T)) = t-T$, where $t \geq T$. Hence for any $t \geq T$, by triangle inequality,  we have the following
$$d_{X}(\Psi(P_n)(t)),\Psi(Q_n)(t))\hspace*{88mm} $$
$$\leq d_{X}( \Psi(P_n)(t),  \Psi(P_n)(T)) + d_{X}(\Psi(P_n)(T),\Psi(Q_n)(T))  + d_{X}(\Psi(Q_n)(T), \Psi(Q_n)(t)) $$
$$\leq \frac{\epsilon}{4\beta} + 2(t-T) \leq \frac{\epsilon}{4} + 2(t-T)\hspace*{68mm}$$

On the other hand, for any $0\leq t \leq T$, by Lemma \ref{control}, we have
$$d_{X}(\Psi(P_n)(t)),\Psi(Q_n)(t)) \leq \beta d_{X}(\Psi(P_n)(T)),\Psi(Q_n)(T)) <  \frac{\epsilon}{4}$$
In particular,
$$d_X(P_n,Q_n) = d_{X}(\Psi(P_n)(0)),\Psi(Q_n)(0)) < \frac{\epsilon}{4}$$
And for $t \leq 0$, $\Psi(P_n)(t) = \Psi(P_n)(0) = P_n$, $\Psi(Q_n)(t) = \Psi(Q_n)(0) = Q_n$, hence
$$d_{X}(\Psi(P_n)(t)),\Psi(Q_n)(t)) =  d_{X}(P_n,Q_n) < \frac{\epsilon}{4\beta}.$$

 Therefore, for any $n>\bar{N}$
$$d_{HFS(X)}(\Psi(P_n),\Psi(Q_n)) = \int_{\mathbb{R}} \frac{d_X(\Psi(P_n)(t)),\Psi(Q_n)(t)))}{2e^{|t|}} dt \hspace*{30mm} $$
$$ = \int_{(-\infty,0]} \frac{d_X(\Psi(P_n)(t)),\Psi(Q_n)(t)))}{2e^{|t|}} dt ~~+~~ \int_{[0,T]} \frac{d_X(\Psi(P_n)(t)),\Psi(Q_n)(t)))}{2e^{|t|}} dt$$

$$\hspace*{20mm}+ \int_{[T, \infty)} \frac{d_X(\Psi(P_n)(t)),\Psi(Q_n)(t)))}{2e^{|t|}} dt$$

$$ \leq \int_{(-\infty,0]} \frac{\frac{\epsilon}{4}}{2e^{|t|}} dt \hspace*{2mm}+ \hspace*{2mm}\int_{[0,T]} \frac{\frac{\epsilon}{4}}{2e^{|t|}} dt  \hspace*{2mm}+\hspace*{2mm} \int_{[T,\infty)} \frac{\frac{\epsilon}{4} + 2(t-T)}{2e^{|t|}} dt \hspace*{20mm}$$
$$\leq \frac{\epsilon}{4} +  e^{-T} \leq \frac{\epsilon}{4} + \frac{\epsilon}{4} = \frac{\epsilon}{2}\hspace*{75mm}$$
Hence we proved the Lemma.
\qed

Because of the properties proved in Lemma \ref{l3}, Theorem 1.4 in \cite{BLR} yields  a long thin cover for $HFS(X)$; i.e. the following result holds

\begin{prop}\label{ltc}
 There exists a natural number $N$, depending only on $k_\Gamma,~d_{HFS(X)}$ and the action of $\Gamma$ on an arbitrary neighborhood of ${HFS(X)}^{\mathbb{R}}$ such that for every $\lambda > 0$ there is an $\mathcal{VC}yc$-cover $\mathcal{U}$ of $HFS(X)$ with the following properties:
 \begin{itemize}
  \item[(i)] dim $\mathcal{U} \leq N$;
 \item [(ii)] For every $x \in HFS(X)$ there exists $U_x \in \mathcal{U}$ such that
 $$ \Phi_{[-\lambda, \lambda]}(x) := \{ \Phi_{\tau}(x) ~|~ \tau \in [-\lambda, \lambda]  \}   \subseteq U_x; $$

  \item [(iii)] ${\Gamma} \setminus {\mathcal{U}}$ is finite.
 \end{itemize}
where $\mathcal{VC}yc$ denote the collections of virtually cyclic subgroups of a group.
 \end{prop}



Recall that the dimension of a cover $\mathcal{U}$ is defined to be the greatest $N$ such that there exists $N+1$ elements in $\mathcal{U}$ with nonempty intersection.  In general, for a collection of subgroups $\mathcal{F}$, we define a $\mathcal{F}$-cover as following.

\begin{defn}
Let G be a group and Z be a G-space. Let $\mathcal{F}$ be a collection of subgroups of G. An open cover $\mathcal{U}$ of Z is called an \textbf{$\mathcal{F}$-cover} if the following three conditions are satisfied.\\
(i) For $g\in G$ and $ U \in \mathcal{U}$ we have either $g(U) = U$ or $g(U) \bigcap U = \emptyset$;\\
(ii) For $g\in G$ and $ U \in \mathcal{U}$, we have $g(U) \in \mathcal{U}$;\\
(iii) For $U \in \mathcal{U}$ the subgroup $G_U := \{ g \in G~|~ g(U) = U \}$ is a member of $\mathcal{F}$.

\end{defn}

For a subset A of a metric space Z and $\delta >0$, $A^{\delta}$ denotes the set of all points $z \in Z$ for which $d(z,A) < \delta$. Combining Lemma \ref{inq} and the fact that $\Gamma$ acts cocompactly on ${FS(T_d) \times \mathbb{R}}$ (Lemma \ref{l3}, (i)), Proposition \ref{ltc} can be improved to the following.

\begin{prop}\label{lte}
 There exists a natural number $N$, depending only on $k_\Gamma,d_{HFS(X)}$ and the action of $\Gamma$ on an arbitrary neighborhood of ${HFS(X)}^{\mathbb{R}}$ such that for every $\lambda > 0$ there is a $\mathcal{VC}yc$-cover $\mathcal{U}$ of $HFS(X)$ with the following properties:
 \begin{itemize}
  \item[(i)] dim $\mathcal{U} \leq N$;
 \item [(ii)] There exists a $\delta >0$ depends on $\lambda$ such that for every $x \in HFS(X)$ there exists $U_x \in \mathcal{U}$ such that
 $$ (\Phi_{[-\lambda, \lambda]}(x))^{\delta}   ~\subseteq~  U_x; ~~~~~~~$$

  \item [(iii)] ${\Gamma} \setminus {\mathcal{U}}$ is finite.
 \end{itemize}
 \end{prop}
\Proof A simple modification of the argument in \cite{BLR}, section 1.3, page 1804-1805, yields the result. In their proof, they used a lemma (Lemma 7.2) which will be replaced by Lemma \ref{inq} in our case. \qed

\section{Proof of the main theorem}

\subsection{Some induction}
Let $D$ be a integral domain with quotient field $K$, let $a \in U(D)$, where $U(D)$ denotes the units of $D$. Form the semi-direct product group $D^+\rtimes_a \BZ$, where $D^+$ denotes the additive group of $D$ and $\BZ$ acts on $D^+$ via multiplication by $a$.

\begin{lem}\label{lim}
If FJC is true for $D^+\rtimes_a \BZ$, then it is also true for $K^+\rtimes_a \BZ$, and hence for every subgroup of $K^+\rtimes_a \BZ$.
\end{lem}
\Proof Let $b \in D \setminus 0$, then
$$D[\frac{1}{b}] = \bigcup_{n=1}^{\infty} \frac{1}{b^n}D^+ = \lim_{n\rightarrow \infty} \frac{1}{b^n} D^+ ,$$
and
$$ ({D[\frac{1}{b}]})^+\rtimes_a \BZ = \lim_{n\rightarrow \infty} {\frac{1}{b^n}D^+}\rtimes_a \BZ $$
Since ${\frac{1}{b^n}D^+}\rtimes_a \BZ \cong D^+\rtimes_a \BZ$, FJC is true for ${\frac{1}{b^n}D^+}\rtimes_a \BZ$ and hence for the direct limit group $ {D[\frac{1}{b}]}^+\rtimes_a \BZ$  by Proposition \ref{qut}(3). On the other hand
$$K^+\rtimes_a \BZ = \lim_{b \in D\setminus 0}  {D[\frac{1}{b}]}^+\rtimes_a \BZ,$$
Hence FJC is true for $K^+\rtimes_a \BZ$. And by \ref{qut}(1), any subgroup of $K^+\rtimes_a \BZ$ also satisfies FJC.
\qed

Let $m(x)$ be a monic polynomial in $\BZ[x]$ which is irreducible in $\BQ[x]$. Let $K$ be the number field determined by $m(x)$; i.e.
$$K = \BQ[x]/m(x)\BQ[x],$$
and let $\mathcal{O}_K$ denote its ring of integers.

The rest of the section is devoted to prove the following result.

\begin{thm}\label{main}
The FJC is true for the group $\Gamma = {\mathcal{O}_K [\frac{1}{x}]}^+\rtimes_x \BZ$, and hence true for  $K^+\rtimes_x \BZ$ by Lemma \ref{lim}.
\end{thm}

\begin{cor}\label{two}
Let $M \in M_n (\BZ)$ have determinant $d>1$. Assume that its characteristic polynomial $m(x)$ is irreducible in $\BQ[x]$, then $ \Gamma = (\BZ[\frac{1}{d}])^n \rtimes_M \BZ$ satisfies FJC. Moreover, $\BQ^n \rtimes_M \BZ$ satisfies FJC.
\end{cor}

\Proof It suffices because of Theorem \ref{main} to embed $\Gamma$ in $K^+\rtimes_x \BZ$. And we now proceed to construct this embedding. Let $e_1,e_2, \cdots, e_n$ be a basis for $\mathcal{O}_K^+$. (Recall that $\mathcal{O}_K^+$ is a free abelian group of rank $n$). Note that $x \in \mathcal{O}_K^+$ since $m(x) = 0$ and $m(x)$ is a monic polynomial with integral coefficients. Hence multiplication by $x$ induces an endomorphism  $f$ of $\mathcal{O}_K^+$ and $f$ determines a matrix $\bar{M} \in M_n(\BZ)$ using basis $\{e_1,e_2, \cdots, e_n\}$. Furthermore the minimum polynomial of $\bar{M}$ is also $m(x)$. So $M$ and $\bar{M}$ have the same minimum and characteristic polynomial; namely $m(x)$ which is $Q[x]$ irreducible. Consequently $M$ and $\bar{M}$ are conjugate via a matrix $T \in GL_n(\BQ)$. Let $\bar{e}_1,\bar{e}_2,\cdots, \bar{e}_n$ be a new basis for $K = \mathcal{O}_K \otimes \BQ$ determined by $T$ and the basis $e_1,e_2, \cdots, e_n$. Then multiplication by $x$ on $K$ determines the matrix $M$ in terms of $\bar{e}_1,\bar{e}_2,\cdots, \bar{e}_n$. We can now define the desired group homomorphism $F: (\BZ[\frac{1}{d}])^n \rtimes_M \rightarrow K^+\rtimes_x \BZ$  by
$$(b,k)  \rightarrow (\bar{b},k),$$
where $(b,k) \in \BZ[\frac{1}{d}]^n \rtimes_M \BZ$, $b = (b_1,b_2,\cdots,b_n)$ with $b_i \in \BZ[\frac{1}{d}]$ and $\bar{b} = b_1\bar{e}_1 +b_2\bar{e}_2 + \cdots b_n\bar{e}_n$. This homomorphism is clearly monic. Since this extends to an embedding of $\BQ^n \rtimes_M \BZ$ into $K^+ \rtimes_x \BZ$, $\BQ^n \rtimes_M \BZ$ also satisfies FJC. \qed

\subsection{Hyper-elementary subgroups of $(\BZ/s\BZ)^n \rtimes_{M_s} \BZ/r\BZ$} \label{gts}
Recall a hyper-elementary group is an extension of a $p$-group by a cyclic group of order $n$ such that $(n,p) = 1$ (compare Definition \ref{he}).

\begin{prop}\label{hgp}
Given $M$ a $n\times n$ matrix with integer coefficients, $det(M) = d > 1$, then we have the following:\\
Given any positive integer $N$, there are positive integers $s$ and $r$ satisfying\\
\textbf{(i)} $s \equiv ~1 ~mod~ d$; \\
 \textbf{(ii)} The order of $GL_n(\BZ/s\BZ)$ divides $r$. In particular we can consider the group
 ${(\BZ / s\BZ)}^n \rtimes_{M_s} \BZ/r\BZ$, where $M_s$ is the reduction of $M$ modulo $s$. Note that $M_s \in GL_n(\BZ/s\BZ)$ since $det(M) = d$ is coprime to $s$.\\
 \textbf{(iii)} If $H$ is a hyper-elementary subgroup of $(\BZ/s\BZ)^n \rtimes_{M_s} \BZ/r\BZ$, then at least one of the following two statements is true:\\
 \hspace*{10mm}(a) $[\BZ/r\BZ, pr(H)] \geq N$, where $pr$ the projection map from $(\BZ/s\BZ)^n \rtimes_{M_s} \BZ/r\BZ$ to $\BZ/r\BZ$.\\
\hspace*{10mm}(b)  There exists a natural number $k$ satisfying \\
           \hspace*{40mm} $k$ divides $s$;\\
           \hspace*{40mm} $k \geq N$;\\
            \hspace*{40mm} $H \bigcap (\BZ/s\BZ)^n \rtimes \{0\} \subset k (\BZ/s\BZ)^n $;

\end{prop}
\Proof The result here is just saying $M$  is hyperbolic good in the sense of Bartels-Farrell-L\"uck in \cite{BFL}, Definition 3.11, by taking $o = d$ which  ensures that the reduction of $M$, $M_s \in GL_n(\BZ/s\BZ)$. The only thing is here $det(M) = d >1$, fortunately the prove of Lemma 3.18 to Lemma 3.21 there still applies directly here. \qed

\subsection{Proof of Theorem \ref{main}}

The remaining proof here for Theorem \ref{main} is exactly parallel to the proof of  the FJC for the solvable Baumslag-Solitar groups in \cite{FW}, Section 5. Hence we will only give a sketch here.

\textbf{Claim} The group ${\mathcal{O}_K [\frac{1}{x}]}^+\rtimes_a \BZ$ is a Farrell-Hsiang group (compare Definition \ref{fhm}) with respect to the family of virtually abelian subgroups.

If we can prove the Claim, by Theorem \ref{fht}, then we proved ${\mathcal{O}_K [\frac{1}{x}]}^+\rtimes_x \BZ$ satisfies the FJC with respect to the family of abelian subgroups. Since FJC is known for abelian groups, by the transitivity principle (Lemma \ref{trp}), we have proven our Theorem \ref{main}. As noted in the proof of Corollary \ref{two}, ${\mathcal{O}}_K$ is a free abelian group of rank $n$, we will fix a basis for it, denote the action of $x$ under this basis by M, which is an integer matrix with determinant $d >1$. Let $M_s$ denote the reduced matrix of $M$ modulo $s$.

\textbf{Sketch of the proof of the Claim} Given any $m > 0$, we will choose the finite group $F_m$ to be $(\BZ/s\BZ)^n \rtimes_{M_s} \BZ/r\BZ$, where $r$ and $s$ depends on $m$, $r$ divides the order of $GL_n(\BZ/s\BZ)$. We need to construct a surjective group homomorphism
$$\alpha_m: {\mathcal{O}_K [\frac{1}{x}]}^+\rtimes_x \BZ \rightarrow (\BZ/s\BZ)^n \rtimes_{M_s} \BZ/r\BZ$$
Note that ${\mathcal{O}_K [\frac{1}{x}]}^+$ is isomorphic to $(\BZ[\frac{1}{d}])^n$, and $s$ is coprime to $d$, $s{\mathcal{O}_K [\frac{1}{x}]}^+ \rtimes \{0\}$ is a normal subgroup of ${\mathcal{O}_K [\frac{1}{x}]}^+\rtimes_x \BZ$ with quotient isomorphic to $(\BZ/s\BZ)^n \rtimes_{M_s} \BZ$. On the other hand $M_s \in GL_n(\BZ/s\BZ)$ while $r$ divides the order of $GL_n(\BZ/s\BZ)$. Hence we can further map $(\BZ/s\BZ)^n \rtimes_{M_s} \BZ$ to $(\BZ/s\BZ)^n \rtimes_{M_s} \BZ/r\BZ$.

 Now by Proposition \ref{hgp} the hyper-elementary subgroup of $F_m$ will be divided into two cases  and $N$ will be determined later, but at least bigger than $4m^2$.

 Case (a), $[\BZ/r\BZ, pr(H)] \geq N > 4m^2$, where $pr$ is the projection map from $(\BZ/s\BZ)^n \rtimes_{M_s} \BZ/r\BZ$ to $\BZ/r\BZ$. In this case we will choose $E_H$ to be the real line $\BR$ with large cells. And the map $f_H$ from ${\mathcal{O}_K [\frac{1}{x}]}^+\rtimes_a \BZ$ to $\BR$ is just the projection map to the second factor as $\BZ$ naturally embeds in $\BR$. For more details, compare the proof in \cite{FW}, Section 5, Case (1).

Case (b), There exists a natural number $k > N$ such that, $H \bigcap (\BZ/s\BZ)^n \rtimes \{0\} \subset k (\BZ/s\BZ)^n $ and $k$ divides $s$. We will use the following diagram to produce $E_H$ in this case.

$\hspace*{10mm}\begin{CD}
\Gamma  @> {\eta} >> T_d \times \mathbb{R}^n @>{F_{k}^{-1}}>> T_d \times \mathbb{R}^n @> {\Psi_\tau}>>  HFS(T_d \times\mathbb{R}^n) \\
\end{CD}$ \\*
~\\*
where $F_{k}: T_d  \times \mathbb{R}^n \rightarrow T_d \times \mathbb{R}^n$ by \\*
$\hspace*{20mm} F_{k} (z,w) = \left(
\begin{array}{cc}
    k & 0 \\
   0 & 1 \\
    \end{array}\right) (z,w)$, for $(z,w) \in T_d\times \mathbb{R}^n$\\*
and $\eta$ is just the embedding by picking up base point $(P_0,0)$, i.e. $\eta(g) = g(P_0,0)$. The action of $\left(
\begin{array}{cc}
    k & 0 \\
   0 & 1 \\
    \end{array}\right)$ on $(z,w) \in T_d\times \mathbb{R}^n$ is the diagonal action explained in section \ref{model}. In particular $F_{k}^{-1}$ will shrink  $T_d  \times \mathbb{R}^n $ in the $\BR^n$ direction. $\Psi_\tau$ is defined in Remark \ref{remfl}. By Proposition \ref{lte}, we have a long thin cover on $HFS(T_d \times\mathbb{R}^n)$, using $\Psi_\tau$  and ${F_{k}^{-1}}$ we can pull back the cover to $T_d \times \BR^n$. $E_H$ will be the nerve of this cover, there is a canonical map $\bar{f}_H: T_d \times \BR^n \rightarrow E_H$, we define $f_H = \bar{f}_H \circ \eta$. Now by some careful choices of $k$ and $\tau$, one can prove that $f_H$ satisfies the inequality we need in the definition of Farrell-Hsiang group (Definition \ref{fhm}). For more details, compare the proof in \cite{FW}, Section 5, Case (2).

~\\
F. Thomas Farrell\\
DEPARTMENT OF MATHEMATICS, SUNY BINGHAMTON, NY,13902 U.S.A.\\
E-mail address: farrell@math.binghamton.edu\\
~\\
Xiaolei Wu\\
DEPARTMENT OF MATHEMATICS, SUNY BINGHAMTON, NY,13902 U.S.A.\\
E-mail address: xwu@math.binghamton.edu\\
Current address: Freie Universit\"at Berlin, Institut f\"ur Mathematik, Arnimallee 7, 14195 Berlin Germany

\end{document}